\documentclass[leqno]{macrorendm}

\usepackage{xcolor}
\usepackage{subfigure}
\usepackage{booktabs}

\volumeyear{xxxx}\yearnumber{x}\volumenumber{xx}

\def\RR{{\mathbb{R}}}
\def\ZZ{{\mathbb{Z}}}
\def \d {{\bf d}}
\def \p {{\bf p}}
\def \t {{\bf t}}
\def \v {{\bf v}}
\def \z {{\bf z}}
\def \x {{\bf x}}
\def \y {{\bf y}}
\def \w {{\bf w}}
\def \Y {{\bf Y}}

\begin{document}

%

%
\articolo[An application of numerical differentiation formulas to fault detection]{An application of numerical differentiation formulas to discontinuity curve detection from irregularly sampled data}{C.~Bracco, O. Davydov, C. Giannelli, A. Sestini\footnotemark[1]}

%

\begin{abstract}
We present a method to detect {discontinuity curves}, usually {called}  \emph{faults}, from a set of scattered data. The scheme first extracts from the data set a subset of points close to the faults. This selection is based on an indicator obtained by using numerical differentiation formulas with irregular centers for gradient approximation, since they can be directly applied to the scattered point cloud without intermediate approximations on a grid. The shape of the {faults} is reconstructed through local computations of regression lines and quadratic least squares approximations. In the final reconstruction stage, a suitable curve interpolation algorithm is applied to the selected set of ordered points {previously associated with each fault.}

\end{abstract}

%
%
%

\section{Introduction}

The detection of discontinuity curves, usually called  \emph{faults},  for bivariate functions is an important issue in several applications, including image processing and geophysics. Consequently, it attracted the attention of several authors, see e.g., \cite{archibald2005,bozzini2013, crampton2005,romros18}. 

The core of any method devoted to the detection of disconituity curves relies on the choice of a suitable operator that can be used to define an indicator for marking the points along the fault. In this paper we consider specific forms of \emph{numerical differentiation formulas} \cite{davydov2018}, that exhibit significant potential in numerical experiments. The adoption of this kind of formulas is motivated by their definition on irregular centers, that can be naturally exploited when working with scattered data of different nature. The formulas are exact for polynomials of a given order and minimize a certain absolute seminorm of the weight vector.  In addition, error bounds depending on a growth function that takes into account the configuration of the centers can be derived.

Our scheme has three key stages that correspond to the detection, identification, and reconstruction phases of the fault curves. In the first stage, we apply the (minimal) numerical differentiation formulas for gradient approximation on a scattered data set to detect a suitable subset of points that lie close to the fault curves. In the second stage of the algorithm, we compute regression lines and quadratic least squares approximations to identify ordered sets of points that describe the shape of the faults. In the final reconstruction stage, any discontinuity curve is obtained as a smooth planar curve, computed through a spline interpolation scheme applied to each set of these ordered points.

The structure of the paper is as follows. Section~\ref{sec:mndf} briefly reviews minimal numerical differentiation formulas, while the fault detection and reconstruction algorithms are presented in Section~\ref{sec:det}. Numerical examples are then presented in Section~\ref{sec:exm}, and Section~\ref{sec:clo} concludes the paper with some final comments and remarks.

\section{Minimal numerical differentiation formulas}\label{sec:mndf}
Referring to \cite{davydov2018} for the details and focusing on the bivariate case which is here of interest, in this section we briefly introduce numerical differentiation
formulas, which are the mathematical tool used for our direct fault detection approach  from scattered data. Note however that  they have a fundamental interest in the
context of meshless numerical solution of PDEs, see references and discussion in \cite{davydov2018}.  

Numerical differentiation formulas are a generalization of finite differences on grids, aimed to approximate the value assumed at a certain site $\z$ from a non vanishing linear differential operator $D$ of order $k$ applied to a sufficiently regular function $f,$ 
\begin{equation} \label{diffop} Df(\z) := 
\sum_{\alpha \in \ZZ_+^2, |\alpha| \le k} c_\alpha(\z)\partial^\alpha f(\z),\quad\text{with}\quad
\partial^\alpha f:=\frac{\partial^\alpha f}{\partial x_1^{\alpha_1}\partial x_2^{\alpha_2}}\,.  
\end{equation}
Such approximation is defined just as a linear combination of the values of $f$ on a finite set $X_\z := \{ \x_i, i=1,\ldots,N_\z\} \subset \RR^2$ of scattered data sites,
\begin{equation}\label{diff1}
\hat Df(\z):=\sum_{i=1}^{N_\z} w_i f(\x_i)\,.
\end{equation}
The meaningfulness of the formula defined in (\ref{diff1}) is ensured by requiring its exactness for bivariate polynomials up to a certain order $q$ (i.e.\ total degree at
most $q-1$).
Actually, as expected thinking to what happens in the context of ordinary differential equations,  Theorem \ref{theoNDF}  reported below and proved in \cite{davydov2018}
underlines that it is necessary to require $q>k$ in order to ensure that the error $\vert \hat Df(\z) - Df(\z) \vert$ tends to zero for sufficiently smooth
$f$ when $h_{\z,X_\z} \rightarrow 0, $ where
$h_{\z,X_\z}$ measures the size of the neighbourhood of $\z$ containing $X_\z$,
\begin{equation} \label{hdef} h_{\z,X_\z} := \max_{1 \le i \le N_\z} \|\x_i - \z\|_2 \,.
\end{equation}

Let $C^{r,\gamma}(Q)$, $r=0,1\ldots$, $\gamma \in (0\,,\,1]$, 
be the H\"{o}lder space consisting of all $r$ times continuously differentiable functions $f$ in $Q$  
with $ \vert \partial^\alpha f \vert_{0,\gamma,Q} < + \infty$, ${\alpha} \in \ZZ^2_+$, $|\alpha| = r$, where
\begin{equation} \label{holdnorm} \vert f \vert_{0,\gamma,Q} \,:=\, \sup_{\x, \y \in Q,\, \x \ne \y} \frac{|f(\x)-f(\y)|}{\|\x - \y\|_2^\gamma}\,. \end{equation}
Furthermore, let $[\x,\y]$ denote the segment joining $\x$ to $\y.$ 

\begin{theorem} \label{theoNDF}
Let $f \in C^{r,\gamma}(\Omega_\z), r=k,\ldots,q-1,$  where  $\Omega_\z$ denotes a domain containing the set
\begin{equation} \label{Sz}   S_\z :=   \displaystyle{\bigcup_{i=1}^{N_\z} [\z\,,\,\x_i] }\,.
\end{equation}
If $D$ is a linear differential operator of order $k$   and $\hat Df(\z)$ is a corresponding numerical differentiation formula with polynomial exactness of order $q>k,$ then
\begin{equation} \label{errorformula} \vert \hat Df(\z) - Df(\z) \vert \le \sigma(\z, X_\z,\w)\, h_{\z,X_\z}^{r+\gamma-k}  \vert f \vert_{r,\gamma,\Omega_\z}\,, \end{equation}
where  $\sigma$ is defined as follows,
\begin{equation} \label{sigmadef}
\sigma(\z, X_\z, \w):= h_{\z,X_\z}^{k-r-\gamma} \sum_{j=1}^{N_\z}\vert w_j \vert \|\x_j-\z\|_2^{r+\gamma},
\end{equation}
and 
$$ \vert f \vert_{r,\gamma,Q} \,:=\,  \frac{1}{(\gamma+1)\cdots (\gamma+r)}\, 
\left( \sum_{\vert \alpha \vert = r}  {r \choose{\alpha}} \vert \partial^\alpha f \vert^2_{0,\gamma,\Omega}  \right)^{1/2} 
\quad\text{if}\quad r\ge1.$$
\end{theorem}

Note that the factor  $\sigma$ defined in (\ref{sigmadef}) is independent of the sample size $h_{z,X}$ \cite{davydov2018}, while it takes into account the cardinality and the shape of $X_\z.$ Besides that, clearly it depends on the specific weight choice.

 \begin{remark}
We observe that, even assuming  $N_\z \ge q(q+1)/2$ (i.e. greater or equal  to the dimension of the space of bivariate
polynomials of order $q$), it may happen that  polynomially exact formulas of order $q>k$ do not exist for some special
geometry of $X_\z$.  This case is surely avoided  if $X_\z$ is {\it unisolvent} for such polynomial space;  if for example 
$N_\z \ge q(q+1)$, then the nonexistence is highly unlikely unless $X_\z$ is a subsample of an unusually irregular distributed data.
\end{remark}

Now in general the $q$--polynomially exact formula  $\hat Df(\z)$ is not unique.  Thus, taking into account that the factor $\sigma$ defined in (\ref{sigmadef}) and appearing in the error bound in (\ref{errorformula}) depends on the specific choice of weights in (\ref{diff1}), it descends that a good strategy for determining the extra degrees of freedom is the first one adopted in \cite{davydov2018} which consists in minimizing the $\ell_1$--weighted (semi)norm of the weights vector $ {\bf w} $ which is defined as $ \vert {\bf w} \vert_{1,r+\gamma} :=  \sum_{i=1}^{N_\z}  \vert w_i \vert~  \Vert \x_i - \z \Vert_2^{r+\gamma}.$
 Conversely, here we choose to minimize the following $\ell_2$--weighted (semi)norm of $\w,$ considering a different proposal also considered in \cite{davydov2018}
\begin{equation}\label{wnorm2}
\vert {\bf w} \vert_{2,r+\gamma} :=\sqrt{\sum_{i=1}^{N_\z}  w_i^2~ \Vert { \x_i}-{\z}\Vert_2^{2(r+\gamma)}}\,.
\end{equation}
 Our choice is justified by the fact that minimizing the $\ell_2$ (semi)-norm is computationally easier and it produces
similar results. Both these formulas are called in the following \emph{minimal numerical differentiation formulas (MNDFs)}.

It is important to note that in both cases MNDFs are \emph{scalable} for homogeneous differential operators \cite{davydov2018} (that is $c_\alpha = 0$ for
all $\alpha$ with $|\alpha| < k$ in \eqref{diffop}), and they are \emph{translation invariant} for operators with constant coefficients.
Therefore for such operators
if we define the following local sample
\begin{equation}\label{Yset} {\bf Y}^h_{\v}(\z) := \{\y_i = \z +\v + h \,(\x_i - \z)\,, \quad i=1,\ldots,N_\z\} \,,
\quad h>0,
\end{equation}
and build an MNDF for the approximation of ${D} f(\z+\v),$  
 
then the corresponding weights $w_i^h$, $i=1,\ldots,N_\z$, of the minimal differentiation formula  do not depend on  
the translation vector $\v$ and they are such that
\begin{equation}
\label{hdip}
w_i^h = w_i /h^k\,, \quad i=1,\ldots,N_\z\,,
\end{equation} 
where $w_i$ are the weights for the case $\v=0$ and $h=1$, {and $k$ is the order of the operator}.

\section{Fault detection and reconstruction} \label{detection}
\label{sec:det}

In this section we introduce the MNDF based indicator we adopt for ordinary fault detection which directly acts on the given
set of scattered data $X \subset \Omega \subset \RR^2$ with associated function values $f(\x)$ for each
$\x\in X$. Analogous indicators related to higher order differential operators
and of interest  for the detection of higher order irregularities of the function $f$ are under development.  

Let us assume that $f \in C^1(\Omega_R),$ with $\Omega_R := \displaystyle{\Omega \setminus \bigcup_{i=1}^M {\cal C}_i}.$ 
Here ${\cal C}_i$ denotes a regular planar curve contained in the domain $\Omega$, without self-intersections and not intersecting anyone of the other
curves. Each ${\cal C}_i$ is   an \emph{ordinary fault} to be detected, that is  we assume that
$$ 
\lim_{\theta \rightarrow 0^+} f(\z+\theta \d_z^L) =: f_z^L \ne f_z^R := \lim_{\theta \rightarrow 0^+} f(\z+\theta \d_z^R)\,, 
\qquad \forall \z \in \mathring{{\cal C}}_i\,,$$
where, denoting by $\t_z$ the unit tangent to ${\cal C}_i$ in $\z,$ $\d_z^L$ and $\d_z^R$ are arbitrary vectors such that 
$\t_z\times \d_z^L<0$ and  $\t_z \times \d_z^R > 0.$ Note that the value of $ \lim_{\theta \rightarrow 0^+} f(\z \pm \theta
\t_z)$ can be either $f_z^L$ or $f_z^R,$ depending on the value of the signed curvature of ${\cal C}_i$ at $\z.$

Denoting by $X_i$ a local subsample of $X$ to be associated with $\x_i \in X,$  
$$X_i:=\{\x_j:\,j\in J_i\}\subset X\,,$$
the \emph{fault indicator} is defined as 
\begin{equation} \label{indic}
I_{\nabla}(\x_i,X_i):=\frac{\Vert \nabla^{i} f (\x_i)\Vert_2}{ \Big\Vert\sum_{ j \in J_i }\vert\w_j\vert \Vert \x_j - \x_i \Vert_2 \Big\Vert_2}\,, 
\end{equation}
relying on the $\ell_2$ MNDF (vector) formula $\nabla^{i} f (\x_i)$ with polynomial exactness $q=2$ approximating $\nabla f(\x_i),$
\begin{equation} \label{indicnum} \nabla^{i} f (\x_i) := \sum_{j \in J_i }\w_j f(\x_j)\,, \end{equation}
where $\w_j$ is a a vector weight, with $\w_j := (w_{j,1}\,,\,w_{j,2})$ and $\vert \w_j \vert := (\vert w_{j,1}\vert\,,\,\vert w_{j,2})\vert).$
Note that the term defining the denominator on the right of (\ref{indic}) does not depend on $f$ and  can be considered 
a scaling factor aimed at reducing the influence of the geometry of $X_i.$  

Concerning this indicator, we observe that it is possible to prove the following two main points \cite{newpaper}. First, the following inequality holds true, 
\begin{equation} 
\label{bound} 
I_{\nabla}(\x_i,X_i) \le \vert f \vert_{0,1,\Omega_i} \,, \qquad \mbox{ if } S_{\x_i} \subset \Omega_i , \mbox{ and } \bar
\Omega_i \subset \Omega_R\,,
\end{equation}
where we note that $\vert f \vert_{0,1,\Omega_i}$ (see the definition in (\ref{holdnorm})) is a measure of the variation of $f$ in a neighborhood $\Omega_i$ of $\x_i.$ 
In particular, $\vert f \vert_{0,1,\Omega_i}\le \max_{\x\in\bar\Omega_i}\|\nabla f(\x)\|_2$ when $\Omega_i$ is convex, and hence 
$I_{\nabla}(\x_i,X_i)$ cannot be too large if $f$ is smooth in $\bar\Omega_i$ with its gradient within
reasonable bounds.  

Furthermore, we can state an asymptotic result related to any point $\x_i \in X.$ Denoting by $\v_i$  any
vector  such that for a nonnegative $\mu_i$ the point ${\bf p}_i  = \x_i - \mu_i \v_i$ belongs to a fault ${\cal C}_k,$ we
need two assumptions on the geometry of $X_i.$ Denoting by $\t_i$ the unit vector tangent to ${\cal C}_k$ at $\p_i,$ first
we require  that there are at least two indices $j_1, j_2 \in J_i$ such that $\left( ( \x_{j_1}-\x_i) \times \t_i \right) \left(
(\x_{j_2}-\x_i) \times \t_i \right)<0.$ Second we need to know that the sum of a subset (suitably defined considering the
geometry of $X_i$) of  the MNDF weights $\w_j$  does not vanish. Then it can be proved that, when $h$ tends to zero, if $\mu_h
\in (0\,,\mu_i)$ is chosen to ensure that $\lim_{h \rightarrow 0} (\mu_i-\mu_h)/h = 0,$ the indicator $I_{\nabla}(\x_i-\mu_h
\v_i , \Y_{-\mu_h \v_i}^h(\x_i))$  tends to infinity,
where $\Y_{-\mu_h \v_i}^h(\x_i)$ is defined according to (\ref{Yset}), {and the formula \eqref{hdip} with $k=1$ applies because the gradient is a 
homogeneous differential operator of order 1 with constant coefficients.}

Therefore, we expect that the indicator value at points close to the fault is much bigger than at points far from it.
Considering this analysis, we say that $\x_i$ belongs to a fault area if $I_{\nabla}(\x_i) > \theta,$ where, clearly, 
the choice of a suitable positive value for the threshold $\theta$ is crucial for the success of the method. 
We can then define the set $F\subseteq X$ of points close to the fault(s) as
\begin{equation}\label{F}
F:=\{\x_i\in X:\, I_{\nabla}(\x_i, X_i) > \theta\}.
\end{equation}
Note that $F$ is a cloud of points concentrated near the faults, from which the shapes of
the faults can be extracted. In
particular, we can use a common approach based on the computation of local least squares approximations (see
\cite{in-kwon2000,soblev18}) to ``narrow'' the point cloud $F$, that is, to obtain a new set $\tilde F$ of points
almost aligned along the faults, from which ordered subsets corresponding to the curves $\mathcal{C}_k$, $k=1,\ldots,M$,
can be extracted using local linear regression lines, and finally approximated, for example with a suitable parametric spline curves, to reconstruct the corresponding fault.

\begin{remark}As well as for other approaches, if there are areas in the considered domain where two faults are very close but the local density of $X$ is not sufficiently high, the presented method cannot detect two separate cloud of points corresponding to the faults. Note however that in this case our method requires only a local increase of the data density to be effective again.
\end{remark}


\section{Examples}
\label{sec:exm}
\begin{figure}[h!]\begin{center}
\includegraphics[trim=0mm 10mm 0mm 10mm,clip,scale=0.50]{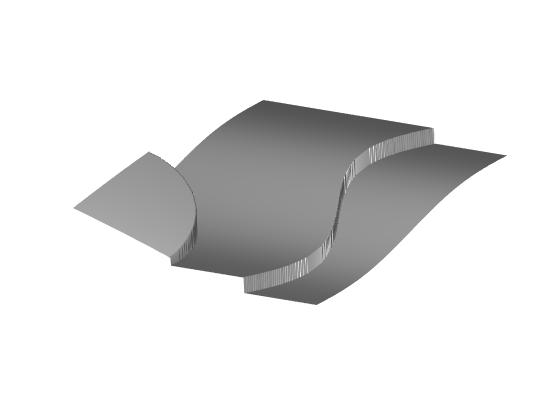}
\caption{Plot of the test function $f$.}
\label{fig1}
\end{center}
\end{figure}

\begin{figure}[h!]\begin{center}
\includegraphics[trim=15mm 5mm 15mm 10mm,clip,scale=0.30]{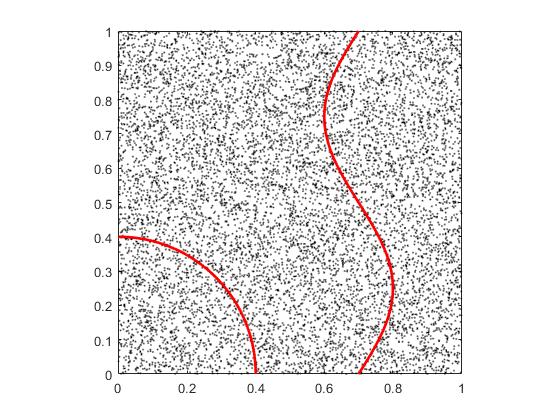}
\includegraphics[trim=15mm 5mm 15mm 10mm,clip,scale=0.30]{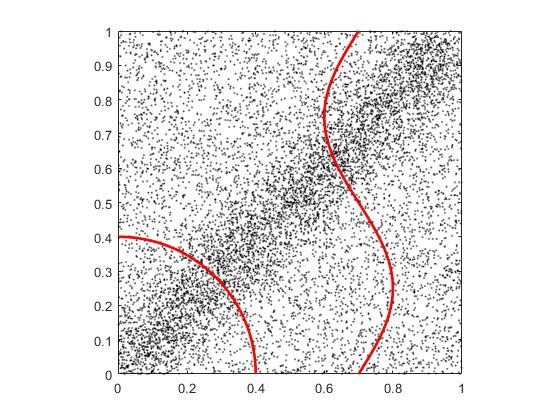}
\\
\includegraphics[trim=15mm 5mm 15mm 10mm,clip,scale=0.30]{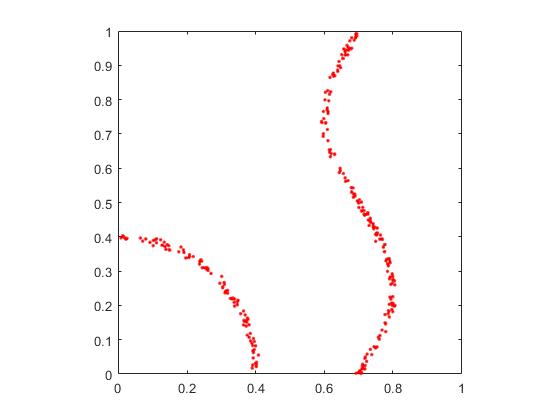}
\includegraphics[trim=15mm 5mm 15mm 10mm,clip,scale=0.30]{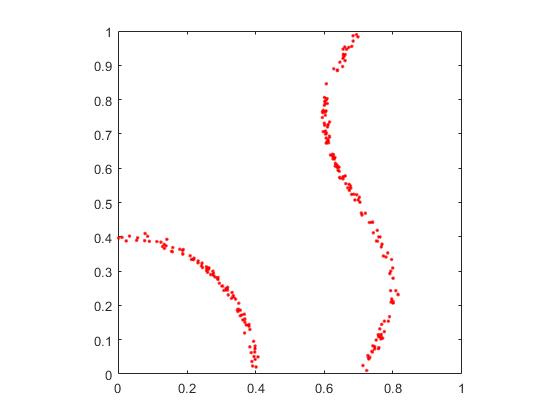}
\\
\includegraphics[trim=15mm 5mm 15mm 10mm,clip,scale=0.30]{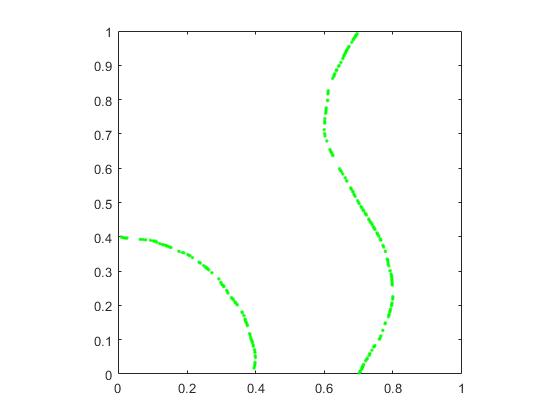}
\includegraphics[trim=15mm 5mm 15mm 10mm,clip,scale=0.30]{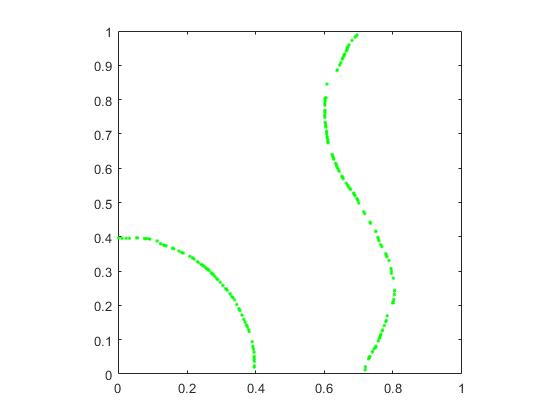}
\\
\includegraphics[trim=15mm 5mm 15mm 10mm,clip,scale=0.30]{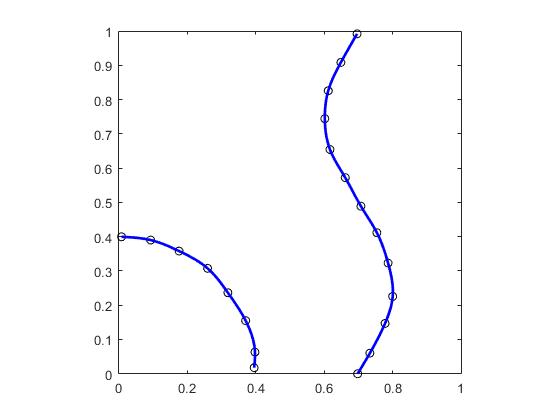}
\includegraphics[trim=15mm 5mm 15mm 10mm,clip,scale=0.30]{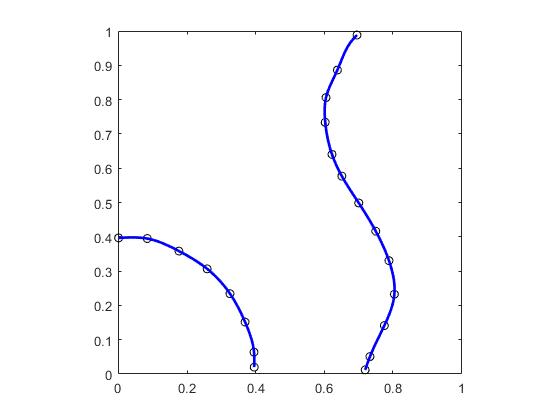}
\caption{Application of the method on a set of  $10000$ uniformly distributed random points (left) and of $9684$ scattered points with varying density (right). 
First row: data sets and exact faults;  second row: detected point sets; third row: narrowed point sets; fourth row: reconstructed faults.
}
 \label{fig2}
\end{center}
\end{figure}

In order to illustrate the presented method, we show two examples of detection and reconstruction of ordinary faults from scattered data corresponding to the test function (see Figure \ref{fig1}) defined on the domain $[0,1]^2$
\begin{equation*}
f(x,y)=\begin{cases} \sqrt{4-x^2-y^2}-2\sqrt{0.99}+0.1 & \mbox{if } x^2+y^2 \le 0.16, \\
  x-0.4-0.1\sin(2\pi y) & \mbox{if } x^2+y^2 > 0.16 \mbox{ and } x \le 0.7 + 0.1\sin(2\pi y),\\
	x-0.4-0.1\sin(2\pi y)-0.2 & \mbox{if } x > 0.7 + 0.1\sin(2\pi y).
\end{cases} 
\end{equation*}
For both examples, {as proposed in Section \ref{detection}, the order of exactness $q$ in the MNDFs used for the fault indicator is set to $2$ and,} for each $\x_i\in X,$ the local sample $X_i$ is chosen as the set of the $q(q+1) = 6$ points of
$X$ closest to $\x_i.$ {Furthermore, the value of the threshold $\theta$ is always set to $0.8$.} \\
 In the first example $X$ is the set of 10000 uniformly distributed random points, see Figure \ref{fig2} (left column, first row). We first determine the subset 
$F\subset X$ of points close to the faults according to \eqref{F}, see Figure \ref{fig2} (left column, second row).
The set $F$ is subsequently \lq\lq narrowed\rq\rq{} to have the points aligned along curves, see Figure \ref{fig2} (left column, third row), and
then ordered subsets of these points are interpolated (in this case by a $C^2$ cubic spline) to
reconstruct the fault curves, see Figure \ref{fig2} (left column, fourth row). Example 2 demonstrates that the method is able 
to handle data of variable density, see the plots on the right column of Figure \ref{fig2}. 
 
In order to assess the accuracy of the results, we also consider the approximated Hausdorff distance between the exact and the reconstructed faults, as well as the maximum distance of the points belonging to the narrowed sets of points from the exact faults. Let $\chi_i$ and $s_i$ be suitable discretizations (composed of $m=500$ points) of the $i$-th exact fault ${\cal C}_i$ and of the computed related approximation, respectively.  Furthermore, let $\tilde F_i$ denote the subset of  narrowed points almost aligned along ${\cal C}_i.$ Then, relating to the considered test function, in Table  \ref{table:tab1} for each of the two faults we report the approximated Hausdorff distance $d_H(s_i,\chi_i)$ and the quantity $d_P(\tilde F_i,\chi_i),$ where
$$
\begin{array}{ll}
d_H(s_i,\chi_i) &:=\, \max\{\displaystyle{\max_{\x \in s_i}\min_{\y \in \chi_i}}\, \Vert \x -\y \Vert_2\,,\,\max_{\y \in \chi_i}\min_{\x \in s_i}\,  \Vert \x -\y \Vert_2\}\,, \cr
\ & \cr
d_P(\tilde F_i,\chi_i) &:=\, \displaystyle{\max_{\ \x \in \tilde F_i}\min_{\y \in \chi_i} }\,  \Vert \x -\y \Vert_2. \cr
\end{array}
$$
 
\begin{table}[!h]

\centerline{ 
\begin{footnotesize}
\begin{tabular}{lcccc} 
\midrule
{}&\multicolumn{2}{c}{${\cal C}_1=$ sinusoid}&\multicolumn{2}{c}{${\cal C}_2=$ quarter of circle} \\
\midrule
{}&{$d_H(s_1,\chi_1)$}& {$d_P(\tilde F_1,\chi_1)$}&{$d_H(s_2,\chi_2)$}& {$d_P(\tilde F_2,\chi_2)$} \\
\midrule
{Example 1}&   $9.102e-3$	& $5.658e-3$ & $1.795e-2$ & $4.871e-3$ 	\\
{Example 2}&   $2.261e-2$	& $1.038e-2$ & $2.048e-2$ & $5.0713e-3$	\\
\midrule
\end{tabular}
\end{footnotesize}
}
\caption{Hausdorff distance $d_H$ between the exact faults and the reconstructed ones, and maximum distance $d_P$ of the points belonging to the narrowed sets of points from the exact faults.}
\label{table:tab1}
\end{table}


%
%

\section{Conclusion}
\label{sec:clo}

We introduced an indicator based on numerical differentiation formulas for the gradient approximation to
detect discontinuity curves from irregularly sampled data. The scheme can be applied to data sets related
to functions with one or more faults. The shape of the faults is then reconstructed in terms of smooth
interpolating splines.

\section*{Acknowledgements}

C.~Bracco, C.~Giannelli, and A.~Sestini are members of the INdAM Research group GNCS. The INdAM support through GNCS and Finanziamenti Premiali SUNRISE is gratefully acknowledged.

\bigskip

\begin{flushleft}

{\bf AMS Subject Classification: 65D10, 65D25}\\[2ex]

%
Cesare~Bracco, Carlotta Giannelli, Alessandra Sestini\\
Dipartimento di Matematica e Informatica ``U.~Dini'', Universit\`a degli Studi di Firenze\\
Viale Morgagni 67a, 50134 Florence, Italy\\
e-mail: \texttt{cesare.bracco@unifi.it, carlotta.giannelli@unifi.it, alessandra.sestini@unifi.it}\\[2ex]
Oleg Davydov\\
Department of Mathematics, Justus Liebig University Giessen\\
Arndtstrasse 2, 35392 Giessen, Germany\\
e-mail: \texttt{oleg.davydov@math.uni-giessen.de}\\[2ex]

%
\textit{Lavoro pervenuto in redazione il MM.GG.AAAA.}

\end{flushleft}

\normalsize
\label{\thechapter:lastpage}

\end{document}